\documentclass[11pt]{amsart}

\usepackage[T1]{fontenc}
\usepackage[utf8]{inputenc}
\usepackage{lmodern}
\usepackage{latexsym}
\usepackage[margin=4.5cm, bottom=3.7cm]{geometry}
\usepackage{indentfirst}
\usepackage{color}
\usepackage{graphicx}
\usepackage{enumerate}
\usepackage{mathtools}
\usepackage{url}

\usepackage{amsmath, amsthm, amssymb, amsfonts}
\newcommand{\N}{\mathbb{N}}

\newcommand{\R}{\mathbb{R}}
\newcommand{\C}{\mathbb{C}}

\newcommand{\SC}{\mathbb{S}}

\newtheorem{tw}{Theorem}[section]
\newtheorem{df}[tw]{Definition}
\newtheorem{pr}[tw]{Proposition}
\newtheorem{wn}[tw]{Corollary}
\newtheorem{lm}[tw]{Lemma}

\theoremstyle{definition}
\newtheorem{example}[tw]{Example}

\theoremstyle{definition}
\newtheorem{rem}[tw]{Remark}

\title[Applications of the Painlev\'e-Kuratowski convergence\dots]{Applications of the Painlev\'e-Kuratowski convergence: Lipschitz functions with converging Clarke subdifferentials and convergence of sets defined by converging equations}
\author{Daniel Fatu\l a}
\address{Jagiellonian University, Faculty of Mathematics and Computer Science, \L ojasiewicza 6, 30-348 Krak\'ow, Poland}\email{daniel.fatula@doctoral.uj.edu.pl}
\date{April 21st, 2024, revised and extended July 13th, 2025}

\begin{document}

\maketitle

\begin{abstract}
In this note we investigate three kinds of applications of the Painlev\'e-Kuratowski convergence of closed sets in analysis that are motivated also by questions from singularity theory. Firstly, we generalise to Lipschitz functions the classical theorem stating that given a sequence of smooth functions with locally uniformly convergent derivatives, we obtain the local uniform convergence of the functions themselves (provided they were convergent at one point). Next, we turn to the study of the behaviour of the fibres of a given function. We prove some general real counterparts of the Hurwitz theorem from complex analysis (stating that the local uniform convergence of holomorphic functions implies the convergence of their sets of zeroes). From the point of view of singularity theory our two theorems concern the convergence of the sets when their descriptions are convergent. They are also of interest in approximation theory and they give some partial results to the problem of when is the limit of a convergent sequence of real algebraic sets algebraic.
\end{abstract}

\section{Introduction}
The Painlev\'e-Kuratowski convergence (sometimes called also Kuratowski convergence) is a natural notion of convergence of closed sets. It generalises the convergence of compact sets in the Hausdorff metric and, as it happens, can be interpreted as the latter, once we use the one-point compactification. As a matter of fact, in the present article we will mostly consider the Painlev\'e-Kuratowski convergence in the Euclidean space, as we are interested in a particular type of applications in analysis in the broad sense. However, in the appendix, we present also a general result, namely we give a new, short and simple proof of the Zarankiewicz Theorem (cf. \cite{b_K}).

The main concern of this paper is to study three problems from singularity theory as proposed by M. Denkowski. The first one is especially interesting on the grounds of the dynamically developing Lipschitz geometry of singularities, but is definitely likely to find also applications in optimization. The question is to obtain a Lipschitz analogue of the classic theorem stating that a sequence of differentiable functions $ f_{n}: \Omega \to \R^{q} $ defined on an open and connected set $ \Omega \subset \R^{p} $ converges locally uniformly if the sequence of their derivatives $ f'_{n} $ converges locally uniformly and there exists a point $ a_0 $ such that the sequence $ f_{n}(a_{0}) $ has a limit. Of course, since Lipschitz functions are only almost everywhere differentiable, we will have to deal with the multifunctions defined by the Clarke subdifferentials. That is where set convergence fits in the picture. We give a complete solution to the problem in Theorems \ref{lip}, \ref{lip2} and the obvious Corollary \ref{lip3}.

The second problem comes directly from singularity theory and asks about the convergence of the sets described by convergent equations. As is often the case, one works with a sequence of functions $f_n$ converging locally uniformly to a function $f$ (e.g. when studying deformations of a given function $f$) and would like to control the behaviour of the sets given by the equations $f_n=0$. The optimal result one would hope to obtain is the convergence of these sets to the set described by the limiting function $f=0$. Observe that this is precisely the kind of question dealt with in the classical Hurwitz Theorem from complex analysis (see \cite{b_TW} for its general version in terms of the convergence of level sets as proposed also in the present article). We give an answer to the question in Theorem \ref{poziomice} where the interesting feature is that even though $f$ is assumed to be $C^1$ smooth, the approximating functions $f_n$ are only required to be continuous. The theorem deals with scalar functions, since in the real case a finite system of equations can always be replaced by a single one (which, however, is a blow to the differential as we need to take a sum of squares thus making the differential to vanish along the described set). The result allows us to study convergence of sets described by equations e.g. algebraic sets. It should also have some repercussions in approximation theory in particular. In this context we give also in the appendix is the negative solution to a question asked by R. Pierzcha\l a by producing a counter-example. Theorem \ref{poziomice} is completed by Theorem \ref{poziomice2} that somehow corresponds to the classical theorem cited above as a motivation for our Lipschitz result Theorem \ref{lip}. Eventually, as an immediate application of these theorems we prove two propositions, \ref{53} and \ref{54}, concerning the problem of obtaining a real counterpart for the algebraic Bishop Theorem of Tworzewski and Winiarski from \cite{TW} stating that the Kuratowski limit of pure-dimensional complex analytic sets with bounded degrees is again algebraic (the only known result up to now is the Pasternak's unpublished solution of the problem for one-dimensional algebraic sets of $\mathbb{R}^2$, \cite{P}).

\section{Preliminaries}

We know from the Rademacher Theorem that a Lipschitz function is differentiable almost everywhere, so that we can define the \textit{generalized differential} or \textit{subdifferential} of $ f: \Omega \to \R^{q} $ at $ x \in \Omega $, where $\Omega\subset\mathbb{R}^n$ is open, in the way Clarke did it in \cite{b_C}:
$$ \partial f(x) := conv \{ \lim_{\nu \to \infty} f'(x_\nu) | x_\nu \in \Omega \setminus Z, x_\nu\to x, \{f'(x_\nu)\}_{\nu=1}^{\infty} \ \textrm{is convergent}\} $$
where $ Z $ is a set of measure $ 0 $ containing the non-differentiability points of $f$ (the above definition does not depend on the choice of $ Z $, see \cite{Cbook} where a paper of J. Warga from 1981 is cited in this respect). The resulting set is a compact and convex subset of the space of linear maps $L(\R^n,\R^q)$.

We also recall Lebourg's version of the Mean Value Theorem for Lipschitz functions proven in \cite{b_L}. If $ f: \Omega \to \R $ is a Lipschitz function and we have two distinct points $ a,b \in \Omega $ such that $ [a, b] \subset \Omega $ then there exist $ c \in (a, b) $ and $ u \in \partial f(c) $ such that
$$ f(b) - f(a) = u(b - a). $$

It is easy to observe that if $ \Omega $ is connected and $ f: \Omega \to \R $ is a Lipschitz function such that $ f'(x) = 0 $ for almost all $ x \in \Omega $ then $ f $ is constant.

Next, for a given metric space $(X,d)$ and a sequence of sets $A_n\subset X$ we define the \textit{Kuratowski lower and upper} limits setting
\begin{align*}
&\liminf_{n\to+\infty} A_n:=\{x\in X\mid \exists x_n\in A_n (n\in\mathbb{N}) \colon x_n\to x\}\\
&\limsup_{n\to+\infty} A_n:=\{x\in X\mid \exists n_1<n_2<\ldots, \exists x_{n_k}\in A_{n_k} (k\in\mathbb{N})\colon x_{n_k}\to x\}.
\end{align*}
It is easy to see that both sets are closed, remain unchanged when we replace $A_n$ with their closures and that $\liminf_{n\to+\infty} A_n\subset \limsup_{n\to+\infty} A_n$. We say that $A_n$ \textit{converges to $A$ in the Painlev\'e-Kuratowski sense}, if both sets coincide with $A$. We shall write then $A_n \xrightarrow{K} A$. More about this convergence can be found in the appendix.

We write $ F: X \rightrightarrows Y $ to denote a \textit{multifunction} $ F: X \to P(Y)$ ($P(Y)$ being the powerset of $Y$). We define \textit{the graph of} $ F $ as $ \Gamma_{F} := \{ (x, y) \in X \times Y|\ y \in F(x) \}. $

There are two natural types of convergence of a sequence of multifunctions (see e.g. \cite{b_RW}).

\begin{df}
Let $ X, Y $ be metric spaces. We take a sequence of multifunctions $ F_{n}: X \rightrightarrows Y $ and $ F: X \rightrightarrows Y. $ We say that
\begin{itemize}
\item $ F_{n} $ tends to $ F $ pointwise if for every $ x \in X $, $ F_{n}(x) \xrightarrow{K} F(x). $
\item $ F_{n} $ tends to $ F $ graphically if $ \Gamma_{F_n} \xrightarrow{K} \Gamma_{F}. $
\end{itemize}
\end{df}

We shall be using the graphical convergence as it is the natural counterpart of local uniform convergence of continuous functions. Besides, the Clarke subdifferentials have closed graphs, as they are outer semi-continuous.

We define the distance function for points $ a \in \R^n $ to a nonempty set $ A \subset \R^n $ as $ dist(a, A) = dist_{A}(a) := \inf \{ |x -a| : x \in A\}, $ where $ |y| $ is the Euclidean norm of $ y \in \R^n. $ We put $ P_{A}(a) := \{ x \in A: |x-a| = dist(x, A)\}$ for the set of closest points to $a$ in $A$.

For $ a, b \in \R^{n} $ we use the notation $ [a, b] := \{ ta + (1-t)b : 0 \leq t \leq 1\}. $

Finally, $B(a,r)$ and $S(a,r)$ will denote the open ball and the sphere centred at $a$ with radius $r>0$ in $\mathbb{R}^n$.

\section{Lipschitz functions with converging Clarke subdifferentials}

Before stating the main theorem, let us provide some enlightening examples concerning the behaviour of the Clarke subdifferentials of a uniformly convergent sequence of Lipschitz functions (compare also Section 7.2 in \cite{BDD}).
\begin{example}\label{np}
Consider the Lipschitz function $f(x)=|x|$ on $[-1,1]$. 
Then 
$$
\partial f(x)=\begin{cases}
    -1,& x\in (-1,0),\\
    [-1,1], &x=0,\\
    1, &x\in (0,1).\end{cases}
$$
Now, let $\tilde{f}_n(x)=\frac{1}{n}f(nx)$ be defined on $\left[-\frac{1}{n},\frac{1}{n}\right]$ and extended to $\mathbb{R}$ by translation: for $k\in\mathbb{Z}$ and $x\in [\frac{2k-1}{n}, \frac{2k+1}{n}]$ we put $f_n(x)=\tilde{f}_n\left(x-\frac{2k}{n}\right)$. All the functions $f_n$ are Lipschitz with the same constant, they converge uniformly to $f_0\equiv 0$, but $\partial f_n$ tends graphically to $\mathbb{R}\times [-1,1]$ which is different from $\partial f_0$ at any point. Note that $\partial f_n$ does not converge pointwise. This can be corrected by using dyadic partitioning.

Not only the graphical limit $D$ of $\partial f_n$ does not have to be almost everywhere univalued, it need not have compact values. To see this, consider $g_n(x)=f(nx)$ on  $\left(-\frac{1}{n},\frac{1}{n}\right)$ and extend these functions continuously by the constant $1$ to the whole of $(-1,1)$. Then 
$$
\partial g_n(x)=\begin{cases}
    0,& \frac{1}{n}<|x|<1,\\
    [-n,n], &x=0,\\
    \mathrm{sgn} x\cdot n, &0<|x|<\frac{1}{n},\\
    \left[\frac{1}{2}(\frac{1}{x}-n),\frac{1}{2}(\frac{1}{x}+n)\right],&x=\pm\frac{1}{n},
    \end{cases}
$$
so that the graphical limit is $D(x)=0$ outside $x=0$ while $D(0)=\mathbb{R}$.
\end{example}

The preceding examples account for the assumptions in the following Theorem. Note also that the graphical convergence corresponds in case of continuous functions to their local uniform convergence.

\begin{tw}\label{lip}
Let $ X $ be an open, connected and bounded subset of $ \R^{p} $ and consider a sequence of mappings $ f_{n}: X \to \R $ such that there exists $ L > 0 $ for which $ | f_{n}(x) - f_{n}(y) | \leq L | x - y | $, for all $ n \in \N $ and $ x, y \in X .$ Assume moreover that the sequence $ \partial f_{n} $ tends graphically to a multifunction $ D: X \rightrightarrows L(\R^p,\R) $ and that $ \# D(x) = 1 $ for almost all $ x \in X. $ Finally, assume that there is a point $ \alpha \in X $ such that $ f_{n}(\alpha) $ is convergent. Then the sequence $ f_{n} $ is uniformly convergent to an $L$-Lipschitz function $ f: X \to \R $ and moreover $ D(x) = \{ f'(x) \} $ for almost all $ x \in X. $
\end{tw}

\begin{proof}
We prove the theorem in a few steps. We use the symbol $ Y $ to denote the dual space of $ \R^{p}. $ We consider the standard norm on $ Y $ defined by $ || u || := \sup_{| x | \leq 1} | u(x) | = \max_{| x | = 1} u(x) $ for $ u \in Y. $ 

Note that all the elements of $\partial f_n(x)$ are bounded by $L$ (the bound for the differentials, wherever they exist, follows directly from the Lipschitz condition and the definition of differentiability), whence all the elements of $D(x)$ are bounded by the same constant.

\textbf{Step 1:} Let $ F $ be a closed subset of $ Y. $ We prove that the set
$$ A_{F} := \{ x \in X: D(x) \cap F \ne \emptyset \} $$
is closed in $ X. $

Consider a sequence $ \{ x_{n} \}_{n=1}^{\infty} \subset  A_F $ which tends to $ x \in X $ and take $ u_n \in D(x_{n}) \cap F. $ Since $ || u_n || \leq L, $ by passing to a subsequence we can assume that $ u_n $ tends to some $ u \in F. $ We take a sequence $ (x_{n, k}, u_{n, k}) \in \Gamma_{\partial f_{k}} $ such that $ \lim_{k \to \infty} (x_{n, k}, u_{n, k}) = (x_{n}, u_{n}). $ So there exist $ k_n \in \N $ large enough such that $ || (x_{n, k_n}, u_{n, k_n}) - (x_{n}, u_{n}) || < 2^{-n} $ and we can additionally assume that $ k_{n+1} > k_{n}. $ Then $ (x, u) = \lim_{n \to \infty} (x_{n, k_n}, u_{n, k_n}) \in \limsup_{n \to \infty} \Gamma_{\partial f_{n}} = \Gamma_{D} $ which ends the proof of this step.

\textbf{Step 2:} We show that if $ K \subset X $ is a compact set, $ V \subset Y $ is open and $ D(x) \subset V $ for $ x \in K, $ then there exists $ N \in \N $ such that for all $ n > N $ and $ x \in K $, $ \partial f_{n}(x) \subset V. $

Suppose that there exist $ x_{n} \in K $ and $ u_{n} \in \partial f_{n}(x_{n}) \setminus V \ne \emptyset $, for infinitely many $ n \in \N. $ Since $ || u_n || \leq L, $ $ (x_{n}, u_{n}) \in \Gamma_{f_{n}} $ has a convergent subsequence to $ (x, u) \in \Gamma_{D}. $ But $ x \in K $ and $ u \notin V $ which gives us a contradiction.

\textbf{Step 3:} In this step we assume that $ f_n $ tends pointwise to a function $ f: X \to \R. $ Obviously $ | f(x) - f(y) | \leq L | x - y| $ for $ x, y \in X. $ We take $ a \in X $ and $ u \in Y $ such that $ D(a) = \{ u \} $ and $ f $ is differentiable at $ a. $ We prove that $ f'(a) = u. $

We take $ y \in \R^{p} $ such that $ || y || =1 $ and $ (u - f'(a))(y) = || u - f'(a) || $ and we fix any $ \epsilon > 0. $
Since $ a \notin A_{Y \setminus B(u, \epsilon)}, $ from Step 1 we get that there exists $ r_{\epsilon} \in (0, \epsilon) $ such that for all $ x \in \overline{B}(a, r_{\epsilon}) $, $ D(x) \subset  B(u, \epsilon). $ From Step 2 we obtain $ N \in \N $ such that for all $ n > N $ and $ x \in \overline{B}(a, r_{\epsilon}) $, we have $ \partial f_{n}(x) \subset  B(u, \epsilon). $
From Lebourg's Mean Value Theorem we have that for $ n \in \N $, there exist $ c_{n}^{\epsilon} \in \overline{B}(a, r_{\epsilon}) $ and $ u_{n}^{\epsilon} \in \partial f_{n}(c_{n}^{\epsilon}) $ such that
$$ f_{n}( r_{\epsilon}y + a) - f_{n}(a) = u_{n}^{\epsilon}(r_{\epsilon}y). $$

Now, $ \{ u_{n}^{\epsilon} \}_{n=0}^{\infty} $ has a convergent subsequence to some $ u^{\epsilon} \in \overline{B}(u, \epsilon). $ It is easy to see that $ f( r_{\epsilon}y + a) - f(a) = u^{\epsilon}(r_{\epsilon}y) $ and
$$ \frac{ f( r_{\epsilon}y + a) - f(a) - f'(a)(r_{\epsilon}y)}{| r_{\epsilon}y |} = (u^{\epsilon} - f'(a))(y). $$

Finally, since $ u^{\epsilon} \in \overline{B}(u, \epsilon), $ by passing with $ \epsilon $ to $ 0 $ we obtain
$$ 0 = (u - f'(a))(y) = || u - f'(a) ||. $$

\textbf{Step 4:} By the Kirszbraun Theorem we know that that all the functions $ f_n $ can be extended to Lipschitz functions defined on $ \overline{X}$ with the same constant $L$. Therefore, if the sequence $ f_n $ were not uniformly convergent, we can use the Arzel\`a-Ascoli Theorem to obtain two subsequences uniformly convergent to some distinct mappings $ f $ and $ g. $ By Step 3 we have that $ D(x) = \{ f'(x) \} = \{ g'(x) \} $ for almost all $x \in X. $ In particular $ (f - g)' = 0 $ almost everywhere. That gives us that $ f - g $ is constant. Since $ \lim_{n \to \infty} f_{n}(\alpha) = f(\alpha) = g(\alpha) $ we see that $ f = g $ which ends the proof.
\end{proof}

\begin{example}
    In the Theorem above we cannot hope to get $D=\partial f$ everywhere in general. To see this, consider the functions $\tilde{f}_n$ from Example \ref{np} and extend them by 0 to $(-1,1)$. Then they tend to $f_0\equiv 0$, but $D(x)=0$ only for $x\neq 0$, while $D(0)=[-1,1]$.
\end{example}

The last Theorem is also true for functions taking values in $ \R^{q}, $ but first we need to prove a lemma.

Let us denote the natural projections $ \pi_{i}: \R^{q} \ni (x_{1}, \ldots, x_{q}) \longmapsto x_{i} \in \R $ and $ \Pi_i : L(\R^{p}, \R^{q}) \ni u \longmapsto \pi_i \circ u \in L(\R^{p}, \R) $ for $ i = 1, \ldots, q. $

\begin{lm}
Let $ X $ be open in $ \R^{p}, $ $ f: X \to \R^q $ be a Lipschitz function and $ a \in X. $ Then $ \partial (\pi_i \circ f)(a) = \Pi_{i}(\partial f(a)). $
\end{lm}

\begin{proof} Let us denote $ f_i := \pi_i \circ f. $ We take $ Z $ to be a subset of $ X $ of measure $ 0, $ such that $ f $ is differentiable for every $ x \in Z' := X\setminus Z. $ We define the sets $ \widehat{\partial}f(x) := \{ \lim_{y_n \to x} f'(y_n) | Z'\ni y_n\to x\colon  \{f'(y_n)\}_{n=1}^\infty\ \textrm{is convergent}\} $ and $ \widehat{\partial}f_{i}(x) := \{ \lim_{y_n \to x} f'_{i}(y_n) |Z'\ni y_n\to x\colon \{f'_i(y_n)\}_{n=1}^\infty\ \textrm{is convergent} \}. $

It is clear from the definition, that $ \Pi_i (\widehat{\partial}f(a)) \subset \widehat{\partial}f_{i}(a). $ For the reverse inclusion take $ u_{i} \in \widehat{\partial}f_{i}(a) $ and $ \{ y_n \} \subset Z' $ such that $ \lim_{n \to \infty} y_n = a $ and $ \lim_{n \to \infty} f'_{i}(y_n) = u_{i}. $ It is easy to notice that for $ j = 1, \ldots, q $ we have $ || f'_{j}(y_n) || \leq L, $ where $ L $ is the Lipschitz constant for $ f. $ We can take $ \{ y_{k_n} \}, $ such that $ f'_{j}(y_{k_n}) $ tends to some $ u_{j} \in (\R^p)' $ for $ j \in \{1, \ldots, q \} \setminus \{ i \}. $ We obtain, that $ (u_{1}, \ldots u_{q}) \in \widehat{\partial}f(a).  $ So $ u_i \in \Pi_i (\widehat{\partial}f(a)) $ and $ \Pi_i (\widehat{\partial}f(a)) = \widehat{\partial}f_{i}(a). $

$ \Pi_i (\partial f(a)) = \Pi_i (conv \widehat{\partial}f(a)) = conv \Pi_i (\widehat{\partial}f(a)) = \partial f_{i}(a) $ follows from the linearity of $ \Pi_i $ and the Carath\'eodory theorem.
\end{proof}

\begin{tw}\label{lip2}
Let $ X $ be an open, connected and bounded subset of $ \R^{p} $ and consider a sequence of mappings $ f_{n}: X \to \R^{q} $ such that there exists $ L > 0 $ for which $ | f_{n}(x) - f_{n}(y) | \leq L | x - y | $, for all $ n \in \N $ and $ x, y \in X .$ Assume that the sequence $ \partial f_{n} $ tends graphically to a multifunction $ D: X \rightrightarrows L(\R^{p}, \R^{q}) $ such that $ \# D(x) = 1 $ for almost all $ x \in X. $ Finally, assume that there exists a point $ \alpha \in X $ such that $ f_{n}(\alpha) $ is convergent. Then the sequence $ f_{n} $ is uniformly convergent to an $L$-Lipschitz function $ f: X \to \R^q $ and moreover $ D(x) = \{ f'(x) \} $ for almost all $ x \in X. $
\end{tw}

\begin{proof} We fix $ i \in \{1, \ldots, q \} $ and we define $ D_{i}(x) := \Pi_{i}(D(x)). $ It is enough to show that $ \partial (\pi_i \circ f_n) $ tends graphically to $ D_{i}. $

Consider $ (x, \Pi_{i}(u)) \in \Gamma_{D_i}, $ where $ u \in D(x). $ There exists a sequence $ (x_{n}, u_{n}) \in \Gamma_{\partial f_{n}} $ convergent to $ (x, u). $ Then $ \lim_{n \to \infty} (x_{n}, \Pi_{i}(u_{n})) = (x, \Pi_{i}(u)). $ So $ \Gamma_{D_i} \subset \liminf_{n \to \infty} \Gamma_{\partial (\pi_i \circ f_n)}. $

Now we want to prove that $ \limsup_{n \to \infty} \Gamma_{\partial (\pi_i \circ f_n)} \subset \Gamma_{D_i}. $ We consider $ (x, v) \in \limsup_{n \to \infty} \Gamma_{\partial (\pi_i \circ f_n)}. $ We obtain a sequence $ \{ (x_{n}, \Pi_{i}(u_n)) \} $ convergent to $ (x, v), $ such that for infinitely many $ n \in \N $, we have $ u_n \in \partial f_{n}(x_n) $ and $ (x_{n}, \Pi_{i}(u_n)) \in \Gamma_{\partial (\pi_i \circ f_n)}. $ Since $ || u_n || \leq L, $ we can pass to a subsequence $ \{ u_{k_n} \} $ convergent to $ u \in D(x). $ So $ (x, v) = \lim_{n \to \infty} (x_{k_n}, \Pi_{i}(u_{k_n})) \in \Gamma_{D_i}. $
\end{proof}

\begin{rem}
As we have seen in the course of the proof, due to the Kirszbraun Theorem, we may as well assume that all the functions $f_n$ are defined on the compact set $\overline{X}$. Even though, due to the Zarankiewicz Theorem, passing to a subsequence guarantees that the limit $D$ can be extended to $\overline{X}$, we can no longer be sure that $D$ will remain univalued almost everywhere ($\partial X$ may have positive measure). 

However, if we assume in the last Theorem that all the assumptions are satisfied on the compact set $\overline{X}$ with $D$ being defined on it, we may somewhat weaken the assumption that all the functions $f_n$ in the sequence admit the same Lipschitz constant $L$. Instead, we may require only that they be Lipschitz functions and that, in addition, the limiting multifunction $D$ be uniformly bounded in the following sense: there should exist $R>0$ such that $\Gamma_D\subset \overline{X}\times \subset B(0,R)$. This in fact implies that all the functions share a common Lipschitz constant.

    Indeed, all the graphs $\Gamma_{\partial f_n}$ are compact and connected. Then the graphical convergence implies that $\Gamma_{\partial f_n}\overline{X}\times \subset B(0,R)$ for all sufficiently large indices. In particular, at each differentiability point $x$ of such an $f_n$ we have $f'_n(x)\in \partial f_n(x)\subset B(0,R)$. It follows that $f_n$ is locally $R$-Lipschitz, but $\overline{X}$ being compact, we conclude that there is a constant $L>0$ depending on $R$ but not on $n$ and such that all the $f_n$'s are $L$-Lipschitz.
\end{rem}

We can also formulate the last Theorem if $ X $ is not bounded, but we obtain only the local uniform convergence.

\begin{wn}\label{lip3}
Let $ X $ be an open and connected subset of $ \R^{p} $ and consider a sequence of mappings $ f_{n}: X \to \R^{q} $ such that there exists $ L > 0 $ for which $ | f_{n}(x) - f_{n}(y) | \leq L | x - y | $, for all $ n \in \N $ and $ x, y \in X. $ Assume that the sequence $ \partial f_{n} $ tends graphically to a multifunction $ D: X \rightrightarrows L(\R^{p}, \R^{q}) $ such that $ \# D(x) = 1 $ for almost all $ x \in X. $ Finally, assume that there exists a point $ \alpha \in X $ such that $ f_{n}(\alpha) $ is convergent. Then the sequence $ f_{n} $ is locally uniformly convergent to an $L$-Lipschitz function $ f: X \to \R $ and moreover $ D(x) = \{ f'(x) \} $ for almost all $ x \in X. $
\end{wn}

\begin{proof} We define $ A := \{ a \in X: f_n (a) \mbox{ is convergent} \}. $ By the last theorem $ A $ is open and closed in $ X, $ so $ A = X. $ We define $ f := \lim_{n \to \infty} f_{n}. $ It follows from the previous Theorem that $ f_n $ tends uniformly to $ f $ on any open, connected and bounded subset, thus $ f_n$ converges to $f $ locally uniformly.
\end{proof}

\section{Convergence of sets defined by converging equations}

In this part we generalize to the real case the classical Hurwitz theorem which says that given an open subset $ \Omega $ of $ \C $ and a sequence of holomorphic mappings $ f_{n}: \Omega \to \C $ tending locally uniformly to a non-constant holomorphic function $ f: \Omega \to \C $ (the holomorphicity follows from the type of convergence, due to the Weierstrass Theorem), if $ z_0 $ is a zero of order $ k $ of $ f $, then there exists $ U $ a neighbourhood of $ z_0 $ such that for all $ n $ big enough, $ f_{n} $ has exactly $ k $ zeros in $ U $ (counted with multiplicities). A generalization of this theorem to higher dimensions was proven in \cite{b_TW}. Now we present two real counterparts of those results.

\begin{tw}\label{poziomice}
Let $ \Omega $ be an open subset of $ \R^{p}. $ Assume that the sequence of continuous mappings $ f_{n}: \Omega \to \R $ tends locally uniformly to a $ C^1 $ function $ f: \Omega \to \R. $ Then for $ b \in \R $, a regular value of $ f $, we have: $ f_{n}^{-1}(b) \xrightarrow{K} f^{-1}(b) $ in $ \Omega. $
\end{tw}

\begin{proof} Take a compact set $ K \subset \Omega $ such that $ a_{n} \in K \cap f_{n}^{-1}(b) $ for infinitely many $ n \in \N. $ Then choose $ a_{k_n} $ a subsequence convergent to $ a \in K. $ We know that $ f_{k_n}(a_{k_n}) \to f(a), $ but $ f_{k_n}(a_{k_n}) = b, $ so $ f(a) = b. $ That means $ a \in K \cap f^{-1}(b) $ and the latter is nonempty.

We only need to prove that for $ a \in f^{-1}(b) $ and $ U \subset \Omega $ a neighbourhood of $ a $, there exist $ N \in \N $ such that for all $ n > N $, we have: $ U \cap f_{n}^{-1}(b) \ne \emptyset. $

We can assume that $ U $ is connected. Since $ b $ is a regular value of $ f, $ we know that $ f $ does not have a local extremum at $ a. $ So there exist $ c_{1}, c_{2} \in U $ such that $ f(c_{1}) < b < f(c_{2}). $ Then we can find $ N \in \N $ such that for $ n > N $ we have: $ f_{n}(c_{1}) < b < f_{n}(c_{2}). $ From the Darboux property, for $ n > N $, we obtain $ a_n \in U $ such that $ f_{n}(a_{n}) = b. $
\end{proof}

\begin{tw}\label{poziomice2}
Let $ \Omega $ be an open set in $ \R^{p}. $ For $ n \in \N $ we consider $ C^{1} $ functions $ f, f_{n}: \Omega \to \R^{q}, $ where $ q \leq p. $ Assume that $ f'_{n} \to f' $ locally uniformly and $ f_{n} \to f $ pointwise. Then for $ b \in \R^{q} $, a regular value of $ f $, we have: $ f_{n}^{-1}(b) \xrightarrow{K} f^{-1}(b) $ in $ \Omega. $
\end{tw}

\begin{proof} We observe first that $ f_{n} \to f $ locally uniformly.

Next, we take a compact set $ K \subset \Omega $ such that $ a_{n} \in K \cap f_{n}^{-1}(b) $ for infinitely many $ n \in \N. $ Then we choose $ a_{k_n} $, a subsequence convergent to $ a \in K. $ We know that $ f_{k_n}(a_{k_n}) \to f(a), $ but $ f_{k_n}(a_{k_n}) = b, $ so $ f(a) = b. $ That means $ a \in K \cap f^{-1}(b) \ne \emptyset. $

We only need to prove that for $ a \in f^{-1}(b) $ and $ U \subset \Omega $, a neighbourhood of $ a $, there exists $ N \in \N $ such that for all $ n > N $, we have: $ \overline{U} \cap f_{n}^{-1}(b) \ne \emptyset. $

Since $ f'(a) $ is an epimorphism, there exists a linear map $ L: \R^{p} \to \R^{p-q} $ such that $ (f'(a), L) $ is an isomorphism. Then we define $ g(x) := (f(x), L(x)), $ $ g_{n}(x) := (f_{n}(x), L(x)) $ and $ b' := g(a). $ Obviously $ g_{n} \to g $ locally uniformly and $ g'_{n} \to g' $ locally uniformly. Then there exists $ \Omega_{0} \subset \Omega $, a neighbourhood of $ a $ such that $ g|_{\Omega_{0}} $ is a diffeomorphism and there exists $ n_{0} \in \N $ such that for $ n > n_{0} $ and $ x \in \Omega_{0} $, $ g_{n}'(x) $ is an isomorphism. In particular $ g_{n}|_{\Omega_{0}} $ are open for almost all $ n \in \N. $

By way of contradiction we fix $ U \subset \Omega_{0} $ a relatively compact neighbourhood of $ a $ and we assume that $ b' \notin g_{n}(\overline{U}) $ for infinitely many $ n \in \N. $ For any such $ n $ we take $ x_{n} \in \overline{U} $ such that $ g_{n}(x_{n}) $ realizes the distance of $ b'$ to $ g_{n}(\overline{U}) $. We observe that $ x_{n} \notin U, $ because $ g_{n}(U) $ is open. Since $ x_{n} \in \partial U $ it has a subsequence $ x_{k_n} $ convergent to $ x \in \partial U. $ From the definition of $ x_n $ we have $ | g_{k_n}(x_{k_n}) - b' | \leq | g_{k_n}(a) - b' | \to 0. $ But $ g_{n} \to g $ locally uniformly which implies $ g_{k_n}(x_{k_n}) \to g(x),$ so $ g(a) = g(x). $ Since $ a \in U $ and $ x \in \partial U, $ we have a contradiction with the injectivity of $ g|_{\Omega_{0}}. $

Summing up, $ b' \in g_{n}(\overline{U}) $, for almost all $ n \in \N. $ Thus $ b \in f_{n}(\overline{U}) $ for those $ n, $ which ends the proof.
\end{proof}

Let us consider an application of the above theorems. In complex geometry it has been known for a long time that the Kuratowski limit of a sequence of pure dimensional complex algebraic sets with uniformly bounded projective degrees is algebraic -- this result of Tworzewski and Winiarski from \cite{TW} may also be seen as an application of the Bishop Theorem on limits of pure dimensional complex analytic sets with locally uniformly bounded volume. The question of under which conditions is a limit of sequence of algebraic sets an algebraic set is definitely very natural. 

The first and apparently only attempt to consider this problem in the real case was made in the master memoir \cite{P}, but has remained unpublished to this day. This is unfortunate as the memoir contains a remarkable result and some very sound observations. In short, firstly, we should expect the limit to be semi-algebraic rather than algebraic (e.g. a sequence of ellipses may converge to a segment); secondly, there is clearly the problem of how to define a degree of a real algebraic set (passing through a complexification seems too complicated, counting the number of points of intersection with a generic affine subspace of codimension equal to the pure dimension of the set, which would correspond to the projective degree from the complex case, does not work, as there is a counterexample in \cite{P}). 

Note that an algebraic subset of $\mathbb{R}^p$ is always described by one equation $P=0$, since any finite system of equations is equivalent to the equation for the sum of their squares. Consider as the degree of an algebraic set $A\subsetneq\mathbb{R}^p$ the number $$\deg A:=\min\{\deg P\mid P\in \mathbb{R}[x_1,\dots, x_n]\colon P^{-1}(0)=A\}.$$ Under the assumption that these degrees are uniformly bounded, Pasternak showed in \cite{P} that a convergent sequence of algebraic one-dimensional subsets of the real plane is semi-algebraic. We can add to this partial result two corollaries from our results.

\begin{pr}\label{53}
    Let $A_n\subset \mathbb{R}^p$ be a sequence of algebraic sets with uniformly bounded degrees, converging to a nonempty set $A$. Then there exist polynomials $P_n\in\mathbb{R}[x_1,\dots, x_p]$ of minimal degree such that $P_n^{-1}(0)=A_n$ and $P_n$ has a subsequence converging locally uniformly to a nonzero polynomial $P$ such that $P^{-1}(0)\supset A$. Moreover, if $0$ is a regular value of $P$, then $P^{-1}(0)=A$, so that $A$ is algebraic.
\end{pr}
\begin{proof}
Let $d \geq \deg A_n$ and $ P_n(x)=\sum_{k \in \N^p, |k| \leq d} a_{k, n}x_1^{k_1}x_2^{k_2}\ldots x_p^{k_p}, $ where for each $n\in \N$ $\max_{k \in \N^p, |k| \leq d}|a_{k, n}| = 1$. After taking a subsequence we can assume that $\lim_{n \to \infty} a_{k,n} = a_k.$ We define $ P(x)=\sum_{k \in \N^p, |k| \leq d} a_{k}x_1^{k_1}\ldots x_p^{k_p}. $ Obviously $P_n$ after passing to a subsequence tend locally uniformly to $P$ and $\max_{k \in \N^p, |k| \leq d}|a_{k}| = 1$, so $P \ne 0$.

To prove $P^{-1}(0)\supset A$ we take $x \in A$. Then we consider a sequence $x_n \in A_n$ such that $\lim_{n \to \infty}x_n = x$. $P_n(x_n)=0$, so $P(x)=0$. The last part is a consequence of the theorem \ref{poziomice}.
\end{proof}

An algebraic subset of $\mathbb{R}^p$ is always described by one equation $P=0$, since any finite system of equations is equivalent to the equation for the sum of their squares, but for this equation 0 is not a regular value, so it can be useful to describe algebraic sets by more than one polynomial equation. It is easy to show the following:

\begin{pr}\label{54}
     Let $A_n\subset \mathbb{R}^p$ be sequence of algebraic sets, converging to a nonempty set $A$. Let $ q\leq p $ and $P_{n,1},\dots, P_{n,q}\in\mathbb{R}[x_1,\dots, x_p]$ be the polynomials that provide the equations for $A_n$. Assume that the degrees of the describing polynomials have a uniform bound. Then for any $n \in \N$ and $i=1, \ldots q$ there exist $a_{n,i} > 0$ such that $a_{n, i}P_{n, i}$ has a subsequence that converges locally uniformly to nonzero polynomials $P_1,\dots, P_q$ such that $P^{-1}(0)\supset A$, where $P= (P_1, \ldots, P_q)$. Moreover, if $0$ is a regular value of $P$, then $P^{-1}(0)=A$, so that $A$ is algebraic.
\end{pr}

\section{Appendix}
In this appendix we gather several additional results that are related to the main topic of this article and seem to be either not well known or maybe new, and of some interest on their own.

\subsection{The Zarankiewicz Theorem}
We present a new, short and simple proof of the Zarankiewicz theorem about the sequential compactness of the Painlev\'e-Kuratowski convergence.

\begin{tw}\label{Zar}
Let $ X $ be a second countable metric space. Consider a sequence of sets $ \{ E_{n} \}_{n=0}^{\infty} \subset P(X). $ Then there exist a subsequence $ E_{k_n} $ convergent to some $ E \in P(X). $
\end{tw}

We will use a set-theoric lemma to prove this theorem.

\begin{df}
For $ A, B \in P(\N) $ we define an equivalence relation: $ A \backsim B $ when there exists a finite set $ S $ such that $ A \cup S = B \cup S. $ We also put $ A^{c} := \N \setminus A. $
\end{df}
\begin{lm}
For any sequence $ \{ A_{n} \}_{n=0}^{\infty} \subset P(\N) $ there exists a sequence $ \{ B_{n} \}_{n=0}^{\infty} \subset P(\N) $ such that $ B_n \backsim A_n $ or $ B_{n} \backsim A_{n}^{c} $ for all $ n \in \N $ and $ \bigcap_{n=0}^{\infty} B_{n} $ is infinite.
\end{lm}

\begin{proof}[Proof of the lemma] Using induction we can take $ C_n = A_n $ or $ C_n = A_{n}^{c} $ such that $ \bigcap_{n=0}^{k} C_{n} $ is infinite for all $ k \in \N. $ Then we consider a sequence $ a_k \in \bigcap_{n=0}^{k} C_{n} $ such that $ a_k > a_{k-1}. $ We define $ B_n := C_n \cup \{ a_{0}, \ldots, a_{n} \}. $ We observe that if $ k \leq n $, then $ a_{k} \in B_{n} $ and if $ k \geq n $, then $ a_k \in C_n \subset B_{n}. $ Thus $ \bigcap_{n=0}^{\infty} B_{n} \supset \{ a_{0}, a_{1}, \ldots \} $ which ends the proof of the Lemma.
\end{proof}

\begin{proof}[Proof of the theorem] We take $ \{ U_{n} \}_{n=0}^{\infty} $ a basis of the topology of $ X $ and we consider the family of sets $ A_{n} := \{ k \in \N| E_{k} \cap U_{n} = \emptyset \} $ for $ n \in \N. $ Let $ \bigcap_{n = 0}^{\infty} B_n = \{ b_{0}, b_{1}, \ldots \}, $ where $ B_n $ are as in the last Lemma. Then we define
$$ E := X \setminus (\bigcup \{ U_{n}| B_{n} \backsim A_{n} \}). $$

We need to prove that $ \limsup_{n \to \infty} E_{b_{n}} \subset E \subset \liminf_{n \to \infty} E_{b_n}. $

Take $ x \in E $ and $ U_{k} $ a neighbourhood of $ x. $ Observe that $ B_{k} \backsim A_{k}^{c}. $ So there exists $ N \in \N $ such that $ \{b_{N}, b_{N+1}, \ldots \} \subset A_{k}^{c}. $ From the definition of $ A_k $ we get that $ E_{b_n} \cap U_{k} \ne \emptyset $, for $ n \geq N. $ In particular $ x \in \liminf_{n \to \infty} E_{b_n}. $

Now, take $ x \notin E, $ then there exists $ U_k $ a neighbourhood of $ x $ such that $ B_k \backsim A_{k}. $ So $ E_{n} \cap U_{k} \ne \emptyset $ only for finitely many $ n \in B_{k}. $ In particular $ x \notin \limsup_{n \to \infty} E_{b_{n}} $ which ends the proof.
\end{proof}

\subsection{The convergence of boundaries implies the convergence of sets in the convex case}

Here we shall give a set-theoretic counterpart of the classical theorem we started with. The boundary of a set can be seen as a kind of `differential' of this set. The natural question about when the convergence of the boundaries implies the convergence of the sets itself leads naturally to restricting the considerations to simply connected sets or, better still, to convex sets. Let us start with some examples.

\begin{example}
    We consider the space $\R^2=\C$ and the sequence of annuli $E_n:=\{z\in\C\mid 1/n\leq|z|\leq 1\}$. Then $\partial E_n\xrightarrow{K} \{|z|=1\}\cup\{0\}$ where the limit set is not the boundary of $\lim E_n=\{|z|\leq 1\}$. 

    A more subtle example is given by the sequence $E_n=nD+(-1)^ni$ where $D$ is the closed unit disc. Then $\partial E_n\xrightarrow{K}\{\mathrm{Im} z=0\}$ while $E_n$ does not converge. It is the fact that the limiting boundary is not compact that destroys the convergence.

    Yet another type of problem is illustrated by an example of $E_n:=\C\setminus\mathrm{int} D$ which yields $\partial E_n=\partial D$, but $E_n$ does not converge to $D$. Note that here the candidate for the limit, i.e. $D$, is compact.

    Finally, let us observe that considering only fat sets (i.e. sets $E$ for which $E=\overline{\mathrm{int} E}$) is a necessary restriction, since we may take a sequence of sets $E_n=\{z\in D\mid (-1)^n\mathrm{Im} z\geq 0\}\cup \partial D$ to get a non-convergent sequence whose boundaries tend to $\partial D\cup [-1,1]$. Of course, $E_n$ are not simply connected.

    We may slightly modify this last example to get simply connected (but not convex) sets $E'_n$ by taking away from the boundary of the set $E_n$ a small arc: $E'_n:=E_n\setminus (\frac{1}{n}\mathrm{int}D+(-1)^{n+1}i)$. Then $\lim\partial E'_n=\lim \partial E_n$, but still $E_n'$ do not converge.

    A tiny modification of $E'_n$ can make these sets become fat. Namely, the arcs sticking out from the half-disc in $E'_n$ may be fattened to take the shape of horns that become thinner and thinner with increasing $n$ so that the limit of the boundaries remains unchanged. A possible formula giving this would be $E''_n:=E'_n\cup(D\setminus \{|z+(-1)^n2i/n|<1-1/n\})$
\end{example}

From the examples above, we see that in the following result the compactness and convexity assumptions are unavoidable and the Theorem seems to be the best one can achieve.

\begin{tw}\label{conv}
Let $ \{ E_n \} $ be a sequence of convex and compact subsets of $ \R^{p}, $ where $ p \geq 2. $ We assume that $ F_n := \partial E_n \xrightarrow{K} F, $ where $ F $ is a compact and nonempty set. Then $ E_n \xrightarrow{K} E = conv F $ and $ \partial E = F. $
\end{tw}

\begin{rem}
    The Theorem is not true for $p=1$. To see this consider the sequence $E_n=[0,n]$ and $E=\{0\}$. But it would be true if we assume that $\bigcup_{n=1}^{\infty} E_n$ is relatively compact. We will see in the next Lemma that this condition is equivalent to the assumption that $F$ is compact a nonempty.
\end{rem}

Before we prove the above theorem we need to state some lemmas.

\begin{lm}
Let $ X $ be a locally compact metric space and let $ \{ E_n \} $ be a sequence of connected subsets of $ X. $ Consider a compact set $ E \in P(X) \setminus \{ \emptyset \} $ and assume that $ E_n \xrightarrow{K} E. $ Then
\begin{enumerate}
\item for every $ \Omega $, neighbourhood of $ E $, there exists $ N \in \N $ such that for all $ n > N $, we have $ E_n \subset \Omega, $
\item there exists $ N \in \N $ such that $ \bigcup_{n > N} E_n $ is relatively compact.
\end{enumerate}
\end{lm}

\begin{proof}
(1) Without loss of generality, we can assume that $ \Omega $ is relatively compact in $ X. $ We take the compact set $ K := \partial \Omega. $ Since $ K \cap E = \emptyset, $ we get $ K \cap E_n = \emptyset $, for all $ n \in \N $ large enough. But that means $ E_n \cap U $ and $ E_n \setminus \overline{U} $ divides $ E_n $ into two open subsets of $ E_{n}. $ By the definition of the Kuratowski lower limit, $ E_n \cap U \ne \emptyset $, for $ n \in \N $ large enough, so $ E_n \subset U $, for $ n \in \N $ large enough.

(2) follows from (1).
\end{proof}

We also recall a well-known lemma (the proof is folklore).

\begin{lm}
Let $ E $ be a compact and convex subset of $ \R^{p}. $ We assume that $ int E \ne \emptyset. $ Then
\begin{enumerate}
\item $ \partial E $ is homeomorphic to $ \SC^{p-1}, $
\item $ E = conv \partial E. $
\end{enumerate}
\end{lm}

Before the proof we can notice the obvious corollary that the boundary of a convex and compact set in $ \R^{p} $ is connected, provided $ p \geq 2. $ Now we are ready to prove Theorem \ref{conv}.

\begin{proof}[Proof of Theorem \ref{conv}]
For a fixed $ \epsilon > 0 $ we consider the open and convex set $ \Omega := \{ x \in \R^{p}: dist(x, E) < \epsilon \}. $ Obviously, $ F \subset \Omega, $ so $ E_n = conv F_{n} \subset \Omega $, for $ n \in \N $ large enough. Then $ \limsup_{n \to \infty} E_n \subset \overline{\Omega}. $ Since $ \epsilon > 0 $ was arbitrary, we obtain $ \limsup_{n \to \infty} E_n \subset E. $
Next, $ \liminf_{n \to \infty} E_n $ is convex as a lower limit of convex sets and $ F \subset \liminf_{n \to \infty} E_{n}. $ So $ E \subset \liminf_{n \to \infty} E_{n}. $ We have proven that $ E_n \xrightarrow{K} E. $

Now we want to prove that $ F \subset \partial E. $ We take $ a \in F $ and we assume that $ a \in int E. $ Then there exist points $ w_{0}, \ldots, w_{p} \in E $ affinely independent and such that $ a $ belongs to interior of $ S = conv \{ w_{0}, \ldots, w_p \}. $ For $ i = 0, \ldots, p $ we take a sequence $ w_{i}^{n} \in E_{n}, $ which tends to $ w_{i}. $ It is easy to observe that $ S_{n} = conv \{ w_{0}^{n}, \ldots, w_{p}^{n} \} \xrightarrow{K} S $ and $ \partial S_{n} \xrightarrow{K} \partial S. $ We fix $ r > 0 $ such that $ \overline{B}(a, r) \subset int S. $ Since the compact set $ \overline{B}(a, r) $ has empty intersection with $ \partial S, $ $ \overline{B}(a, r) \cap \partial S_{n} = \emptyset $, for $ n \in \N $ large enough. So $ \overline{B}(a, r) \subset int S_{n} \subset int E_{n} $, for $ n \in \N $ large enough. But $ B(a, r) \cap \partial E_{n} \ne \emptyset $ for $ n \in \N $ large enough, because $ a \in B(a, r) \cap F \ne \emptyset. $ This contradiction ends our proof of the inclusion $ F \subset \partial E. $

Conversely, assume  $ x \in \partial E \setminus F. $ Then $ dist(x, F) > r $ for some $ r > 0. $ So $ dist(x, F_{n}) > r $, for $ n \in \N $ large enough. If $ x \notin int E_{n} $ for infinitely many $ n \in \N, $ then for those $ n \in \N $ we have $ r < dist(x, \partial E_{n}) = dist(x, E_{n}) \to 0. $ Thus $ x \in int E_{n} $, for $ n \in \N $ big enough. Since $ x \in \partial E, $ then exist a sequence $ x_n \in \R^{p} \setminus E $ which tends to $ x. $ Since $ \{ x_{n} \} \cap E = \emptyset, $ then $ \{ x_{n} \} \cap E_{k} = \emptyset $ for $ k \geq k_{n}. $ We can assume that $ k_n $ tends increasingly to infinity. We take $ z_{n} \in [x, x_{n}] \cap \partial E_{k_n} \ne \emptyset. $ Since $ z_{n} \to x, $ we have $ x \in F, $ which ends the proof.
\end{proof}

\subsection{Remark about the convergence of fibres}

The following question was asked by prof. Rafa{\l} Pierzcha\l a: given an analytic function $f\colon U\to \R$, where $U\subset\R^n$ is a neighbourhood of the origin, satisfying $f(0)=0$, is it true that we can always find a sequence $y_n\to 0$ for which $f^{-1}(y_n)\xrightarrow{K} f^{-1}(0)$? 

As observed in \cite{DD}, the continuity of $f$ already implies the inclusion $\limsup_{y\to 0} f^{-1}(y)\subset f^{-1}(0)$. On the other hand, from \cite{DL} we know that an equivalent condition for $f^{-1}(0)\subset \liminf_{y\to 0} f^{-1}(y)$ to hold in the case of a continuous function is its openness onto its image. Of course, an arbitrary analytic function need not be open onto its image, which makes the question all the more interesting. 

Nevertheless, already in the polynomial case there is a counter-example, namely: $ f(x, y) = x(x - y)^{2}. $

\section{Acknowledgements}
The author is grateful to M. Denkowski for suggesting the problem, several examples concerning it and editorial help.

\section{Additional information}

On behalf of all authors, the corresponding author states that there is no conflict of interest.

\end{document}